\documentclass[reqno,11pt]{amsart}

\usepackage{epsf}
\usepackage{graphics}
\usepackage{graphicx}
\usepackage{amssymb}
\usepackage{amsmath}

\date{}

\theoremstyle{plain}
\newtheorem{theorem}{Theorem}
\newtheorem{corollary}{Corollary}

\newtheorem{rem}{Remark}
\newtheorem{question}{Question}

\theoremstyle{definition}

\theoremstyle{remark}

\newtheorem*{remarks}{Remarks}

\def\N{{\mathbb N}}
\def\Z{{\mathbb Z}}

\def\F{{\mathbb F}}

\title{Monomial reduction of knot polynomials}

\author{Sebastian Baader}

\dedicatory{For Filip, Peter, and Lukas, who all turned 35 recently}

\begin{document}

\begin{abstract}
For all natural numbers $N$ and prime numbers $p$, we find a knot~$K$ whose skein polynomial $P_K(a,z)$ evaluated at $z=N$ has trivial reduction modulo~$p$.
An interesting consequence of our construction is that all polynomials $P_K(a,N)$ (mod~$p$) with bounded $a$-span are realised by knots with bounded braid index. As an application, we classify all polynomials of the form $P_K(a,1)$ (mod~$2$) with $a$-span $\leq 10$.
\end{abstract}

%\dedicatory{}

\maketitle

\section{Introduction}

The HOMFLY polynomial $P_K(a,z) \in \Z[a^{\pm 1},z^{\pm 1}]$ of links~$K$ satisfies the skein relation
$$a^{-1} P_{\widehat{\beta \sigma_i^2}}(a,z)-aP_{\widehat{\beta}}(a,z)=zP_{\widehat{\beta \sigma_i}}(a,z),$$
for all braids $\beta \in B_n$, where $\sigma_i \in B_n$ denotes any standard generator in the braid group $B_n$, and the hat symbol denotes the standard closure of a braid~\cite{HOMFLY}. It is not known whether the property $P_K(a,z)=1$ characterises the trivial knot. The tables of small knots do not even exclude the possibility that $P_K(a,z)$ modulo any prime number detects the trivial knot.

\begin{theorem}
\label{monomial}
For all natural numbers $N$ and prime numbers $p$, there exists a knot $K$ with
$$P_K(a,N)=1 \text{ (mod $p$)}.$$
\end{theorem}

A similar statement is known for the Jones polynomial, but only for finitely many primes~\cite{EF}. The integer substitution $z=N$ transforms the HOMFLY polynomial into a single variable polynomial, whose span provides a lower bound for the minimal braid index~$b(K)$ of knots, by the inequality of Franks-Williams and Morton~\cite{FW,M}:
$$\text{a-span}(P_K(a,z)) \leq 2b(K)-2.$$

\begin{theorem}
\label{span}
For all $N,S \in \N$ and prime numbers $p$, there exists~$T \in \N$, so that the set of polynomials $P_K(a,N)$ (mod~$p$) with $a$-span $\leq S$ is realised by knots~$K$ with $b(K) \leq T$.
\end{theorem}

The above result is of practical use for listing knot polynomials modulo primes, as we will see in the case $p=2$. In particular, we obtain the following concrete result.

\begin{corollary}
\label{span10}
All polynomials of the form $P_K(a,1)$ (mod~$2$) with $a$-span $\leq 10$ are realised by knots with braid index $\leq 7$.
\label{6braids}
\end{corollary}

Thanks to the inequality of Franks-Williams and Morton, the polynomials $P_K(a,1)$ (mod~$2$) with $a$-span $\leq 10$ include the polynomials of all knots with braid index $\leq 6$, hence of all knots with crossing number $\leq 11$.

The main number theoretic input for the first theorem is the fact that every prime number divides a Fibonacci number. The second theorem is more of a coincidental by-product of the construction. Detailed proofs are presented on the next couple of pages.

\section{HOMFLY monomials}

The skein relation for the HOMFLY polynomial allows for an inductive calculation of the polynomials $P_{T(2,n)}(a,z)$ for all torus links $T(2,n)$, defined as closures of the braids $\sigma_1^n \in B_2$ with two strands. The result is a polynomial of $a$-span two. More precisely, for all $k \in \N$, we obtain polynomials of the form
$$P_{T(2,2k)}(a,z)=a_k(z)a^{2k-1}+b_k(z)a^{2k+1},$$
$$P_{T(2,2k+1)}(a,z)=c_k(z)a^{2k}+d_k(z)a^{2k+2},$$
for suitable $a_k,b_k,c_k,d_k \in \Z[z^{\pm 1}]$ with initial values 
$$a_0=z^{-1}, \, b_0=-z^{-1}, \, c_0=1, \, d_0=0.$$
A slight rearrangement of the skein relation,
$$P_{T(2,n+2)}(a,z)=a^2 P_{T(2,n)}(a,z)+azP_{T(2,n+1)}(a,z),$$
yields the following recursive formulas:
$$a_k=a_{k-1}+zc_{k-1},\, b_k=b_{k-1}+zd_{k-1},$$
$$c_k=c_{k-1}+za_k,\, d_k=d_{k-1}+zb_k.$$
Now comes the key observation: after the substitution~$z=N$, and the reduction modulo~$p$, there remain only finitely many possible quadruples
$$(a_k,b_k,c_k,d_k) \in (\Z/\Z_p)^4.$$
As a consequence, there exists $m>0$ with 
$$(a_{k+m},b_{k+m},c_{k+m},d_{k+m})=(a_k,b_k,c_k,d_k) \in (\Z/\Z_p)^4,$$
for all $k \in \N$. In particular, $c_m=1$ and $d_m=0$, thus
$$P_{T(2,2m+1)}(a,N)=a^{2m} \in \Z/\Z_p[a^{\pm 1}].$$
Finally, we observe that the symmetry $$P_{K^*}(a,z)=P_K(a^{-1},z)$$ between the HOMFLY polynomial of a knot~$K$ and its mirror image $K^*$, and the multiplicative rule for connected sums of knots, $P_{K_1 \# K_2}=P_{K_1} P_{K_2}$, imply that the connected sum of knots
$T(2,2m+1)\#T(2,-2m-1)$ has trivial polynomial:
$$P_{T(2,2m+1)\#T(2,-2m-1)}(a,N)=1 \text{ (mod p)}.$$

\begin{remarks} \quad

\noindent
(i) The above argument is inspired by the known fact that every prime number divides a Fibonacci number. In fact, in the special case $N=1$, we obtain
$$P_{T(2,2k+1)}(a,1)=F_{2k+2}a^{2k}+F_{2k}a^{2k+2},$$
where $(F_n)$ denotes the Fibonacci sequence starting at $F_0=0$, $F_1=1$.
The periodicity of the Fibonacci sequence modulo $p$ can be derived from the following explicit result:
$$F_{p-\left(\frac{p}{5}\right)}=0 \text{ (mod p)},$$
where $\left(\frac{p}{5}\right)$ denotes the Legendre symbol~\cite{Le}.

\noindent
(ii) The case $N=0$ (mod $p$) requires a slight adjustment, because of the term $z^{-1}$. An elementary induction shows
$$P_{T(2,2k-1)}(a,0)=ka^{2k-2}-(k-1)a^{2k},$$
for all $k \geq 1$. In particular, $P_{T(2,2p+1)}(a,0)=a^{2p}$ (mod~$p$), from which Theorem~\ref{monomial} follows, as above.
\end{remarks}

As for Theorem~\ref{span}, suppose we are given a polynomial $P_K(a,N) \in \F_p[a^{\pm 1}]$ whose $a$-span is $\leq S$. By the above argument, for all $l \in \N$, we have
$$P_{T(2,2lm+1)}(a,N)=a^{2lm}, \, P_{T(2,-2lm-1)}(a,N)=a^{-2lm}.$$
For a suitable choice of $l$ and sign, we 
can arrange for the minimal $a$-degree of $a^{\pm 2lm}P_K(a,N)$ to be in the range $[0,1,\ldots,2m-1]$. We emphasize that this new polynomial is also the HOMFLY polynomial of a knot $\tilde{K}$, the connected sum of $K$ with a suitable torus knot with two strands. Now there are only finitely many potential polynomials 
$P_{\tilde{K}}(a,N)$ modulo~$p$ with $a$-span $\leq S$ and minimal degree in the range $[0,1,\ldots,2m-1]$. Choose finitely many knots that realise all these polynomials, and let $T-1 \in \N$ be the largest minimal braid index among those. Then the initial polynomial $P_K(a,N) \in \F_p[a^{\pm 1}]$ is realised by a knot with braid index $T$, again by adding a suitable torus knot with two strands to one of these finitely many reference knots with braid index $\leq T-1$.

\section{Small span modulo two reductions}

The goal of this section is to list all knot polynomials $P_K(a,1)$ of $a$-span $\leq 10$, with coefficients modulo two. As in the previous section, the main ingredient is a family of knots with monomial reduction $P_K(a,1)$ (mod $2$). The periodicity of the Fibonacci sequence modulo two being $m=3$, we obtain:
$$P_{T(2,6l+1)}(a,1)=a^{6l}, \, P_{T(2,-6l-1)}(a,1)=a^{-6l},$$
for all $l \in \N$. Therefore, it is enough to classify all knot polynomials $P_K(a,1)$ in the degree ranges $[-4,6]$, $[-2,8]$, $[0,10]$. Indeed, we may discard odd powers of the variable $a$, since the HOMFLY polynomial of knots is in fact a Laurent polynomial in $a^2$ (and a polynomial in $z^2$, for that matter). Moreover, by the mirror symmetry $K \leftrightarrow K^*$, the polynomials in the range $[0,10]$ are in correspondence with the polynomials in the range $[-10,0]$, in turn with the polynomials in the range $[-4,6]$, by multiplying with $a^{6}=P_{T(2,7)}(a,1)$. 
A priori, there are $2^6=64$ potential polynomials $P_K(a,1) \in \F_2[a^{\pm 2}]$ in each degree range $[-4,6]$ and $[-2,8]$. However, the skein relation for the HOMFLY polynomial implies  $P_K(a,a^{-1}-a)=1$, for all knots~$K$. This fact can be enhanced to the following statement~\cite{GGMRS}: the difference between the HOMFLY polynomials of two knots is divisible (in $\Z[a^{\pm 2},z^{\pm 2}]$) by
$$z^2-(a^{-1}-a)^2.$$
After the substitution $z=1$ and reduction modulo two, the latter becomes $a^{-2}+1+a^2$. An elementary argument, based on polynomial division with remainder, shows that every even degree range of span 10 contains only $16$ polynomials $P_K(a,1) \in \F_2[a^{\pm 2}]$ with the additional property that $P_K(a,1)-1$ is divisible by $a^{-2}+1+a^2$. An analogous argument is carried out in~\cite{ABF}, where we classify Jones polynomials modulo two of span eight.

With the help of the knot tables by Rolfsen and knotinfo~\cite{Ro,LM}, we found $16$ polynomials $P_K(a,1) \in \F_2[a^{\pm 2}]$ in each degree range $[-4,6]$ and $[-2,8]$. All the knots can be chosen with braid index $\leq 6$, as seen in Table~1 (where we use the notation of the cited tables, and $O$ stands for the trivial knot). By adding one summand of the form $T(2,\pm (6l+1))$ to all these knots, and by considering all mirror images of the resulting knots, we obtain all polynomials of span $\leq 10$. This implies the statement of Corollary~\ref{span10}. 

\begin{table}
\begin{tabular}{| c |p{9.5cm}|}
\hline
degree range  &  knot polynomials  \\ \hline 
$[-4,6]$  &  $O$ ($1$), $3_1$ ($a^{2}+a^{4}$), $3_1^*$ ($a^{-4}+a^{-2}$), $4_1$($a^{-2}+a^{2}$), $5_1$ ($a^{6}$), $6_1$ ($a^{-2}+1+a^{4}$), $6_1^*$($a^{-4}+1+a^{2}$), $8_3$ ($a^{-4}+a^{-2}+1+a^{2}+a^{4}$), $10_3$ ($a^{-4}+a^{-2}+a^{2}+a^{4}+a^{6}$), $11a103$($a^{-4}+a^{2}+a^{6}$), $11n101$ ($a^{-2}+a^{4}+a^{6}$), $3_1 \# 4_1$ ($1+a^{2}+a^{4}+a^{6}$), $3_1 \# 6_1^*$ ($a^{-2}+1+a^{2}+a^{6}$), $3_1^* \# 11n101$ ($a^{-4}+1+a^{4}+a^{6}$), $4_1 \# 4_1$ ($a^{-4}+a^{4}$), $4_1 \# 6_1$ ($a^{-4}+a^{-2}+1+a^{6}$) \\ \hline
$[-2,8]$  & $O$ ($1$), $3_1$ ($a^{2}+a^{4}$), $4_1$ ($a^{-2}+a^{2}$), $5_1$ ($a^{6}$), $6_1$ ($a^{-2}+1+a^{4}$), $7_4$ ($a^{2}+a^{6}+a^{8}$), $10_1$ ($a^{-2}+1+a^{2}+a^{4}+a^{8}$), $11a121$ ($a^{-2}+a^{2}+a^{4}+a^{6}+a^{8}$), $11n101$ ($a^{-2}+a^{4}+a^{6}$), $11n139$ ($1+a^{2}+a^{8}$), $3_1 \# 3_1$ ($a^{4}+a^{8}$), $3_1 \# 4_1$ ($1+a^{2}+a^{4}+a^{6}$), $3_1 \# 6_1$ ($1+a^{4}+a^{6}+a^{8}$), $3_1 \# 6_1^*$ ($a^{-2}+1+a^{2}+a^{6}$), $3_1 \# 8_3$ ($a^{-2}+a^{8}$), $3_1 \# 4_1 \# 4_1$ ($a^{-2}+1+a^{6}+a^{8}$)\\ \hline
\end{tabular}
\smallskip
\caption{Knot polynomials $P_K(a,1)$ (mod $2$) of span $\leq 10$}
\end{table}

The list of knots obtained in this section suggests the following generalisation.
\begin{question}
Let $p(a) \in \F_2[a^{\pm 2}]$ with $p(a)-1$ divisible by $a^{-2}+1+a^2$. Is there a knot $K \subset S^3$ with $P_K(a,1)=p(a)$ (mod $2$)?
\end{question}

We conclude with a prime reduction of the famous unknot detection question for the HOMFLY polynomial.
\begin{question}
Let $p$ be a prime number. Is there a knot $K \subset S^3$ with
$$P_K(a,z)=1 \text{ (mod $p$)?}$$
\end{question}

%\smallskip
%\noindent
%\texttt{sebastian.baader@unibe.ch}

\end{document}